\documentclass{article}
\usepackage[english]{babel}
\usepackage{amsmath,amssymb,graphicx,latexsym,theorem,bm,mathtools}
\usepackage[citecolor=blue]{hyperref}
\newcommand\blfootnote[1]{%
  \begingroup
  \renewcommand\thefootnote{}\footnote{#1}%
  \addtocounter{footnote}{-1}%
  \endgroup
}
\usepackage{rotating}
\usepackage{float}
\newcommand{\gt}[1]{\bm{#1}}
\renewcommand{\[}{\begin{equation}}
\renewcommand{\]}{\end{equation}}
\newcommand{\R}{\mathbb{R}} 
\newcommand{\tp}[1]{\langle #1 \rangle} 
\newcommand{\bn}[1]{\bar{#1}}
\renewcommand{\[}{\begin{equation}}
\renewcommand{\]}{\end{equation}}
\newcommand{\dd}{\mathrm{d}}
\newcommand{\ten}[1]{\mathbf{#1}} 
\newcommand{\average}[1]{\left\langle#1\right\rangle}

\usepackage{hyperref}
\hypersetup{
    colorlinks=true,
    urlcolor=cyan,
    linkcolor=blue
}

\begin{document}

\title{Complete closed-form solutions to the problem of inextensional bending for surfaces of translation and origami tessellations}

\author{
  Adam Reddy$^1$, Asma Karami$^2$ \& Hussein Nassar$^2$\footnote{Corresponding author: \href{mailto:hnassar@uh.edu}{hnassar@uh.edu}}
}

\date{\small $^1$Department of Mechanical and Aerospace Engineering,
  University of Missouri, Columbia MO 65211\\
 $^2$Department of Mechanical and Aerospace Engineering, University of Houston, Houston TX 77204\\[2ex]%
    \today}
\maketitle

\begin{abstract}
  Plates generally admit six deformation modes: three of which are high in strain energy, stretch the plate's midsurface and are called membrane modes; and three are low-energy, bend the midsurface without stretching it and are called bending modes. For origami tessellations, and other corrugated compliant thin shells, the modes are mixed and it is no longer clear what modes, if any, are low in energy in the sense that they are inextensional. Here, it is shown, by direct construction of closed-form solutions, that when the midsurface is a surface of translation, there exists three infinitesimally inextensional deformation modes that correspond to (1) stretching, with an effective Poisson's effect; (2) bending, with an effective synclastic or anti-clastic effect; and to (3) twisting. The provided expressions are valid irrespective of surface regularity and, in particular, properly handle any creases be them straight or curved. The results provide a powerful benchmark for the validation of numerical methods and further insight into the elastic stiffness of thin corrugated compliant shells.
    \\
    
\end{abstract}

\blfootnote{Work supported by the NSF under CAREER award No. CMMI-2045881}

\section{Introduction}
The isometric deformations of a surface are deflections that preserve distances as measured on the surface. If the surface is the midsurface of a thin shell of thickness $h$ and typical radius of curvature $R$, then an isometric deformation has a strain energy proportional to $(h/R)^3$ all due to bending and with zero membrane contributions. By contrast, deformations that cause stretching have a strain energy dictated by membrane effects and proportional to $h/R$. Then, for small thicknesses ($h/R\ll 1$), isometric deformations dominate the response of a shell, as they are less energetically costly. Now some surfaces do not admit any isometric deformations because of their geometry, topology, or the conditions imposed at boundaries, but other surfaces do. Depending on that, the corresponding shells will be more or less stiff. This suggests that tailoring surfaces with isometric deformations is a purely geometric gateway to designing shells with compliant deformation modes. Compliant shells appear in structural designs~\cite{RADAELLI201778}, modern robots~\cite{Nijssen}, and biological structures~\cite{Flaum} because of a set of desirable attributes that rely on compliance: they can be stowed compactly, deploy into 3D doubly curved geometries, as well as change shape and adapt to changing environments. Foldable origami-like structures in particular have been used to design antennae, satellites, and solar arrays~\cite{WANG2022, Chen2019, Salazar2017}. At smaller scales, compliant or foldable structural elements have been used to achieve exotic mechanical properties in metamaterials~\cite{Meza, Li2017}.

The modeling of thin compliant shells is challenging, as is the underlying geometric problem of finding isometric deformations. Shells that are a discrete assembly of planar facets, as typically found in origami-inspired structures~\cite{Reis}, are exceptional in that they can be treated as a linkage~\cite{FONSECA2022, Seffen2012}. Their kinematics are tractable even if somewhat cumbersome and often employ a mix of vector algebra and spherical trigonometry. Using such techniques, Schenk and Guest~\cite{Schenk2013} and Wei et al.~\cite{Wei2013} were able to determine a fundamental identity between the effective curvature and effective stretching of the Miura origami pattern. Eidini and Paulino~\cite{Eidini2015} and Zhou et al.~\cite{Zhou2016} expanded on these results by examining a general class of patterns that derive from the Miura ori. A broader theory of origami deformations is developed by Brunck et al.~\cite{Brunck2016}.

Beyond these cases, i.e. for shells that are smoothly curved or have curved creases, finding isometric deformations involves finding an appropriate surface parametrization in the sense of differential geometry. For smooth surfaces, many classical results are available for special classes of surfaces, e.g., developable, axisymmetric and Weingarten. More recent contributions provide existence results in the context of curved-crease origami~\cite{Karami2024}, i.e., for piecewise smooth developable surfaces. In particular, Liu et al.~\cite{LIU2024105559} derived a fairly general result which allows one to construct isometric deformations for two surfaces meeting at a crease. More specifically, under mild smoothness assumptions, given a crease and a surface normal defined along the crease, surfaces on both sides of the crease are determined. Demaine et al.~\cite{Demaine} developed several theorems regarding creases and the rulings they produce, including singular cases where kinks form, which they applied to a specific tesselation featuring such kinks. Alternatively, Mundilova~\cite{MUNDILOVA2019} presented an approach that does not rely on examining the Frenet frame and its associated scalar quantities along the crease. Instead, given two distinct developable surfaces, one glues together patches by choosing appropriate boundary curves and rulings from each surface. The caveat to such constructions is that the resulting folding may contain self-intersections or rulings of infinite length depending on the boundary curves and rulings specified. 

Infinitesimal isometric deformations are relaxed isometries in the sense that they preserve distances to leading order in deflections, i.e., for small deflections. They are not as important for morphing and deployment applications, a priori, but are crucial for stiffness computations. For origami tessellations and other periodic surfaces in particular, the infinitesimal isometries identified in the literature are typically ones that maintain periodicity. Indeed, it has been recognized, since the original works of Schenk and Guest~\cite{Schenk2013} and Wei et al.~\cite{Wei2013}, that deflections that represent out-of-plane bending of origami tessellations typically require a non-periodic distribution of fold angles. The logic is that periodic fold angles necessarily lead to cylindrical states~\cite{Tachi2015}. Thus, the Miura-ori, the eggbox pattern and any other origami tessellations that bend into doubly curved forms, be them saddle-like or dome-like, necessarily exhibit non-periodic fold angles, making such states harder to characterize. In the works mentioned above~\cite{Schenk2013, Wei2013}, this was overcome by discretizing quad panels into two triangular facets thereby adding a set of ad-hoc, non-physical, diagonal creases, one per panel. Now Filipov et al.~\cite{Filipov2017} showed, through a series of numerical simulations, that folding about diagonals is indeed appropriate to describe large panel deflections but is not suitable for small panel deflections where other, smoother, kinematics take over. In fact, the analysis affords more refinement: In their study, Filipov et al.~\cite{Filipov2017} simulated a crease pattern with only a few crease lines where it is possible for a load to engage the large deformations of a panel. With large crease patterns, it turns out that the local kinematics remain governed by the infinitesimal deflections of panels even when large deflections are imposed globally through the boundaries~\cite{Nassar2024}. Hence, be it for morphing (i.e., large deflections) or stiffness (i.e., small deflections) computations, the understanding of infinitesimal inextensional kinematics is crucial, at least for surfaces with large repetitive patterns such as origami tessellations.

The infinitesimal isometries of origami tessellations that influence the macroscopic stiffnesses of an equivalent Kirchhoff-Love plate (e.g., tensile, shear, bending and torsion) are necessarily limited in number: there are six of them that are linearly independent at most, three of which correspond to in-plane modes and the other three corresponding to out-of-plane modes~\cite{Sab2020}. Recent theoretical results show in fact that, for simply connected surfaces (i.e., no holes, no handles), there exists exactly three deformation modes, in the sense of Kirchhoff-Love theory of plates, that are infinitesimally isometric~\cite{Nassar2025}. The main contribution of the present paper is to construct closed-form expressions of these three isometries for a broad class of surfaces of translation. It is crucial to highlight that the solutions are valid regardless of whether the surface has flat or curved panels and/or curved or straight crease lines: they simply maintain the regularity of the underlying surface, be it smooth or piecewise smooth. In particular, the solutions do not employ extra crease lines, diagonal or otherwise, and should therefore be more energetically advantageous. Recall that surfaces of translation are swept by the translation of one curve (the profile) along another (the path). They constitute a rich family that includes many well-known origami patterns (e.g., the Miura ori and eggbox patterns) as well as their curved-crease and smooth generalizations. Beyond its theoretical value, the result provides a useful benchmark for the assessment of reduced-order models (e.g., bar and hinge models~\cite{Filipov2017}). The result also simplifies reconstruction procedures of origami tessellations from their limit surfaces (see, e.g., \cite{XU2025106295}).

In what follows, attention is restricted to periodic surfaces for which effective modes of deformation can be easily defined. In the absence of periodicity, the obtained solutions remain valid but their interpretation becomes cumbersome. Each solution is obtained based on an ansatz where variables are separated in a suitable fashion; and each solution corresponds to a mode of deformation with a well-defined effective measure of strain, be it an effective membrane strain or an effective bending strain. This constitutes the complete characterization of (infinitesimally isometric) deformation modes mentioned above. Note that the solutions are unique modulo correctors that are periodic. The choice of corrector does not affect the effective measures of strain and is ultimately fixed by non-geometric considerations, i.e., material and crease properties. In section~2, some preliminary constructions are presented. The main results are derived in section~3, followed by a discussion in section~4 and some concluding remarks.


\section{Theory}
Consider a piecewise smooth surface that describes a thin shell. Piecewise smoothness allows for sudden changes in the tangent plane that model both straight and curved crease lines. Suppose the surface admits a parametrization $\ten r: \R^2 \rightarrow \R^3$ with a metric
\begin{equation}\label{eq:metric}
\ten I  =
   \begin{pmatrix}
        \tp{\ten r_x,\ten r_x} & \tp{\ten r_x,\ten r_y} \\
        \tp{\ten r_x,\ten r_y} &  \tp{\ten r_y,\ten r_y}
   \end{pmatrix}
\end{equation}
where the indices denote partial derivatives of the position vector relative to some curvilinear coordinates $(x,y)$, e.g., $\ten r_x=\partial\ten r/\partial x$. A small deflection given by a displacement $\ten d$ causes the metric to change by amounts given by the infinitesimal strains
\begin{equation}\label{eq:infStrain}
\delta 
\ten I =
   \begin{pmatrix}
	 2\tp{\ten r_x,\ten d_x} &  \tp{\ten r_y,\ten d_x}+\tp{\ten r_x,\ten d_y} \\
	\tp{\ten r_y,\ten d_x}+\tp{\ten r_x,\ten d_y} & 2\tp{\ten r_y,\ten d_y}
   \end{pmatrix}.
\end{equation}
The deformation is isometric, to leading order in displacement amplitude, if $\delta\ten I=\ten 0$, i.e.,
\[\label{eq:infIso}
\begin{split}
    \tp{\ten r_x,\ten d_x}=0,\quad
    \tp{\ten r_y,\ten d_y}=0,\quad
    \tp{\ten r_y,\ten d_x}+\tp{\ten r_x,\ten d_y}=0.
\end{split}
\]

\begin{figure}[ht]
     \centering
     \includegraphics[width=\linewidth]{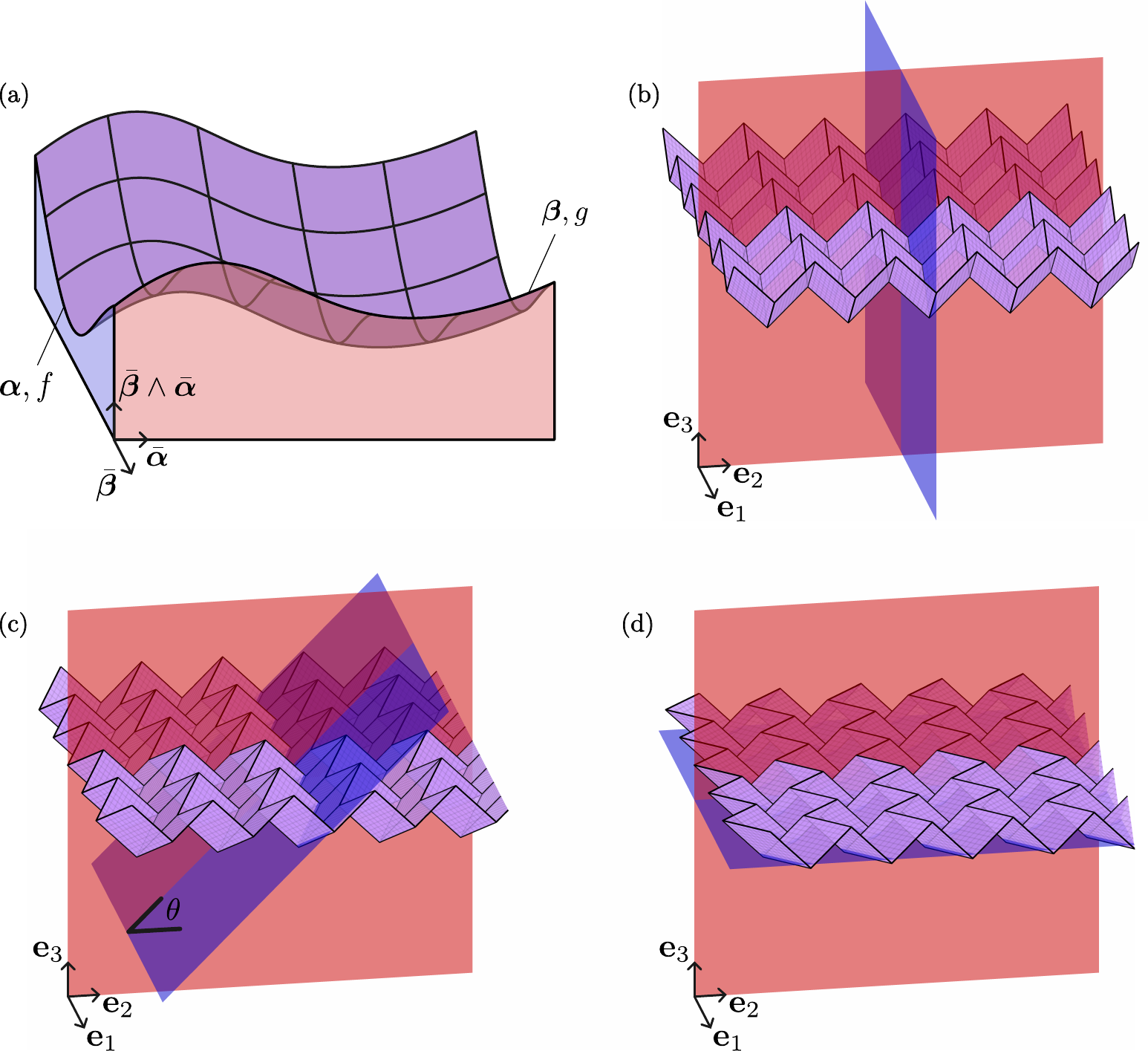}
     \caption{Surfaces of translation, annotated, in (a) the adapted basis; and, (b-d) the canonical basis, with $\theta=0$, $\pi/4$ and $\pi/2$ respectively.}
    \label{fig:SOT}
\end{figure}

Now let $\ten r$ describe a surface of translation spanned by translating one curve~$\gt\alpha$ along another $\gt\beta$; see Figure~\ref{fig:SOT}a. The simplest of such surfaces is of the form
\[
\ten r(x,y) = (x,y,\alpha(x)+\beta(y))
\]
where $(x,y)$ are now cartesian coordinates and $\alpha$ and $\beta$ are two functions whose graphs outline the profiles of the surface in the $x$ and $y$ directions. Here, we afford more generality: $\gt\alpha$ and $\gt\beta$ can feature overhangs that cannot be parametrized using cartesian coordinates. However, we maintain three assumptions:
\begin{enumerate}
    \item each curve belongs to one plane (e.g., helices are excluded);
    \item the tangents to the curves are never in the planes intersection; and,
    \item the planes are perpendicular.
\end{enumerate}
Note that assumptions (1) and (2) are essential. Assumption (3) can be relaxed thanks to the affine invariance of (infinitesimal) isometric deformations but is kept for convenience. Consequently, the surfaces of translation considered in what follows admit parameterizations of the form
\[
    \ten r(x,y) = \gt\alpha(x) + \gt\beta(y),
\]
where the tangents to $\gt\alpha$ and $\gt\beta$ are given by
\begin{equation}\label{eq:surfOfTrans}
\gt\alpha' = \alpha_{1}\bn{\gt\beta}+\alpha_{2} \bn{\gt\beta}\wedge\bn{\gt\alpha}, \quad
\gt\beta' = \beta_{1}\bn{\gt\alpha}+\beta_{2}\bn{\gt\beta}\wedge\bn{\gt\alpha},
\end{equation}
with $\bn{\gt\alpha}$ being the unit normal to the plane of $\gt\alpha$, and similarly for $\bn{\gt\beta}$, so that $[\bn{\gt\beta},\bn{\gt\alpha},\bn{\gt\beta}\wedge\bn{\gt\alpha}]$ is an orthonormal basis as depicted on Figure~\ref{fig:SOT}a. The coefficients $\{\alpha,\beta\}_{\{1,2\}}$ are the coordinates of the tangent vectors in that basis. A prime will generally denote the derivative of a function of a single variable.

The displacement $\ten d$ can conveniently be represented in the same basis:
\[
    \ten d = u\bn{\gt\beta} + v\bn{\gt\alpha} + w \bn{\gt\beta}\wedge\bn{\gt\alpha}.
\]
Then, the isometric constraints~\eqref{eq:infIso} expand into
\begin{equation}\label{eq:isoOfTrans}
\begin{gathered}
\alpha_{1} u_{x} + \alpha_{2}w_{x} = 0,\quad
\beta_{1}v_{y} + \beta_{2}w_{y} = 0,\\
\alpha_{1}u_{y} + \beta_{1}v_{x} + \alpha_{2}w_{y} + \beta_{2}w_{x} = 0,
\end{gathered}
\end{equation}
where again indices $(x,y)$ denote partial derivatives. This is a system of three first-order linear partial differential equations to be solved for the displacement components. Hereafter, we derive three solutions that vaguely correspond to twisting, stretching and bending of the surface $\ten r$. The correspondence becomes rigorous when the profiles $\gt\alpha$ and $\gt\beta$ are periodic.

Results are exemplified for a family of surfaces spanning curved variations on the well-known Miura-ori and eggbox patterns. To define the example surfaces, in reference to the canonical basis $(\ten e_1,\ten e_2,\ten e_3)$ of $\R^3$, let $\gt\alpha$ be the graph of a function $f$ in the plane $(\ten e_1,\ten e_3)$ rotated through an angle $\theta$ about $\ten e_1$, namely
\[
    \gt\alpha(x) = (x,-s f(x), cf(x)),\quad c \equiv \cos\theta,\quad s \equiv \sin\theta.
\]
Let $\gt\beta$ be the graph of a function $g$ in the plane $(\ten e_2,\ten e_3)$, namely
\[
    \gt\beta(y) = (0,y,g(y)).
\]
Then the example surfaces are of the form
\[\label{eq:ExSurf}
    \ten r(x,y) = (x, y - sf(x), g(y) + cf(x)),
\]
and include the particular cases listed in Table~\ref{tab:ExSurf}. See also Figure $\ref{fig:SOT}$. In what follows, it will be useful to note that
\[\label{eq:alphabeta}
\alpha_1 = 1,\quad \alpha_2 = f',\quad \beta_1 = c + sg',\quad \beta_2 = -s+cg'.
\]
Examples are listed in Table~\ref{tab:ExSurf} and reference Figure~\ref{fig:SOT}.
\begin{table}[h!]
    \centering
    \begin{tabular}{c|cccc}
        surface/pattern & $\theta$ & $f$ & $g$ \\ \hline\hline
        Eggbox & $0^\circ$ & zigzag & zigzag & Fig.~\ref{fig:SOT}b\\
        smooth variant & $0^\circ$ & sinusoidal & sinusoidal & Fig.~\ref{fig:SOT}b \\ \hline
         Miura ori & $90^\circ$ & zigzag & zigzag & Fig.~\ref{fig:SOT}d\\
        curved-crease variant & $90^\circ$ & sinusoidal & zigzag & Fig.~\ref{fig:SOT}d\\ \hline
        ``morph'' pattern & $\theta$ & zigzag & zigzag & Fig.~\ref{fig:SOT}c
    \end{tabular}
    \caption{Example surfaces spanned by Equation~\eqref{eq:ExSurf} and depicted on Figure~\ref{fig:SOT}.}
    \label{tab:ExSurf}
\end{table}

\section{Analytical results}
The solutions are derived based on three different schemes of variable separation applied to component $w$, namely
\[
    \begin{array}{l|l|l|l}
    \text{Mode} & \text{Twisting} & \text{Stretching} & \text{Bending} \\
    w(x,y) & X(x)Y(y) & X(x) + Y(y) & X(x)Y(y) + S(x)+T(y)
    \end{array}
\]
where $X$, $Y$, $S$ and $T$ are functions to be determined. Recall that component $w$ is directed along the intersection of the planes of the profiles $\gt\alpha$ and $\gt\beta$.

\subsection{Twisting}
The classical solution of pure torsion for an elastic plate parametrized by $(x,y,0)$ is $(0,0,xy)$. That is, the out-of-plane deflection is linear in the distance traversed along the two axes. In our case, taking these axes to be along $\bn{\gt\beta}$ and $\bn{\gt\alpha}$, the distances traveled are $\int^x\alpha_1(s)\dd s$ and $\int^y\beta_1(s)\dd s$. This leads to the guess $w(x,y)=\int^x\alpha_1(s)\dd s\int^y\beta_1(s)\dd s$ and is of the separated form $X(x)Y(y)$ as promised. It is crucial to note that the components $\{\alpha,\beta\}_{\{1,2\}}$ are discontinuous at crease lines, but their antiderivatives remain continuous. In the following, the notation for antiderivatives, namely $\int^{x,y}\cdot(s)\dd s$, will be shortened into $\int^{x,y}$ or even $\int$, with the integration details being clarified by context.

Substituting $w$ into~\eqref{eq:isoOfTrans} leads to
\[
    u_x = -\alpha_2\int\beta_1, \quad
    v_y = -\beta_2\int\alpha_1,
\]
and, by integration, to
\[
    u = -\int\alpha_2\int\beta_1 + u_o, \quad
    v = -\int\beta_2\int\alpha_1 + v_o,
\]
where $u_o=u_o(y)$ and $v_o=v_o(x)$ are determined by setting the shear component to zero, i.e., thanks to the third equation in~\eqref{eq:isoOfTrans}. Upon rearranging, the latter reads
\[
    \alpha_1\left(\beta_2\int\beta_1-\beta_1\int\beta_2+u_o'\right)+\beta_1\left(\alpha_2\int\alpha_1-\alpha_1\int\alpha_2
    +v_o'\right) = 0.
\]
The separation of variables then implies that there is a constant $\lambda$ such that
\[
u_o' = -\lambda\beta_1+\beta_1\int\beta_2-\beta_2\int\beta_1,\quad
v_o' = \lambda\alpha_1+\alpha_1\int\alpha_2-\alpha_2\int\alpha_1.
\]
The constant is set to zero ($\lambda=0$) as it ends up spanning an in-plane infinitesimal rotation. Integration, and some tidying, lead to the solution
\[
\begin{split}
    u &= -\int\alpha_2\int\beta_1 + \int\left(\beta_1\int\beta_2-\beta_2\int\beta_1\right),\\
    v &= -\int\beta_2\int\alpha_1 + \int\left(\alpha_1\int\alpha_2-\alpha_2\int\alpha_1\right),\\
    w&=\int\alpha_1\int\beta_1.
\end{split}
\]
Note that expressions involving a double integral such as $\iint \alpha_1$ are to be interpreted as a functions such as
\[
\iint\alpha_1: x \mapsto \int^x\dd t\int^t\alpha_1(s)\dd s,
\]
where the lower bound of integration can be chosen arbitrarily, but consistently throughout.

\begin{sidewaysfigure}[!ht]
     \centering
     \includegraphics[width=\linewidth]{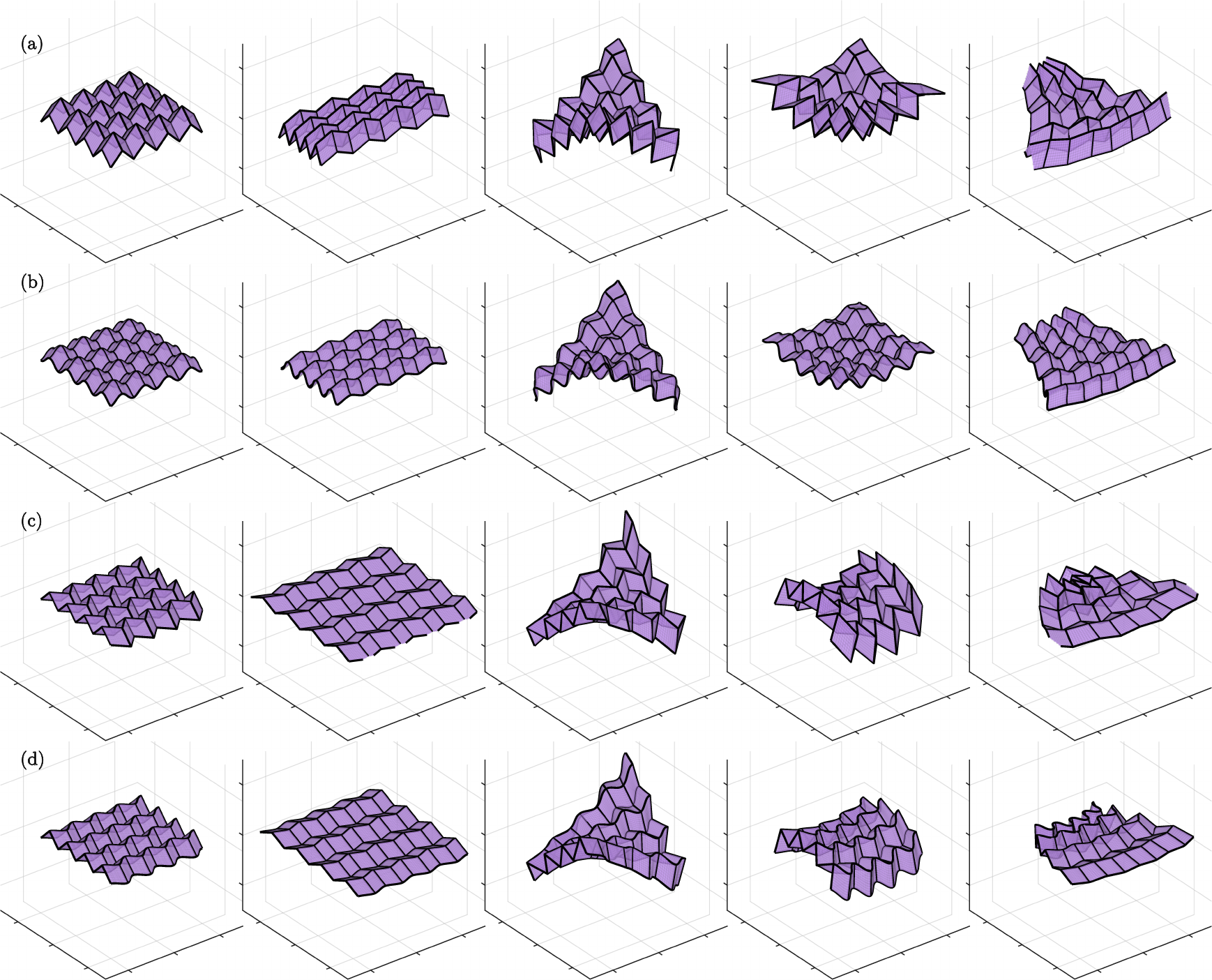}
     \caption{Infinitesimal isometries illustrated for (a) the eggbox pattern, (b) its smoothed version, (c) the Miura ori, (d) its curved-crease version. Modes, per column, left to right are: reference, stretching, twisting, out-of-plane bending and in-plane bending. The Matlab code that produced these figures is available at \url{github.com/nassarh}.}
    \label{fig:deflections}
\end{sidewaysfigure}

For the adopted examples, as described by equation~\eqref{eq:ExSurf}, the deflection is best written in the canonical basis of $\R^3$ where it reads
\[
    \ten d = \begin{pmatrix} 
    -y(g(y)+cf(x)) - sg(y)f(x) + 2G(y) \\
    -x(g(y)+cf(x)) + cF(x) \\
    x(y - sf(x))+2sF(x)
    \end{pmatrix}
\]
where $F$ and $G$ are antiderivatives of $f$ and $g$. Deflected states for selected surfaces (Table~\ref{tab:ExSurf}) are illustrated in Figure~\ref{fig:deflections}. It is remarkable that the vertical deflection, in direction $\ten e_3$, is exactly that of the plane (i.e., $xy$) for $\theta=0^\circ$ irrespective of how the surface is corrugated or creased. In other cases ($\theta\neq 0^\circ$), a correction is needed to ensure the deflection is inextensional.

The above solution is valid whether $f$ and $g$ are periodic or not, but suppose they are. Then, the twisting deformation can be written as
\[
\ten d = \begin{pmatrix}
    -\bar{z}\bar{y} \\ -\bar{z}\bar{x} \\ \bar{x}\bar{y}
\end{pmatrix} + \text{ periodic \& rigid motions}
\]
where $(\bar{x},\bar{y},\bar{z})=(x,y-sf(x),g(y)+cf(x))$ are the cartesian coordinates of the surface. Thus, for a periodic surface, the twisting solution indeed describes the twisting solution in the sense of Kirchhoff-Love plate theory~\cite{Sab2020}.

\clearpage

\clearpage

\subsection{Stretching}
Surfaces of translation, as defined in the present context, can deform isometrically in such a way that they appear to be stretching overall. The unfolding of the Miura ori is a great example of this \emph{effective} stretching. To describe solutions of effective stretching, a variable separation scheme $w(x,y) = X(x)+Y(y)$ is used as announced above. Substituting back into the isometric constraints~\eqref{eq:isoOfTrans} yields
\[
\alpha_1 u_x + \alpha_2 X' = 0, \quad \beta_1 v_y + \beta_2 Y' = 0, \quad
\alpha_1 u_y + \beta_1 v_x + \alpha_2Y' +\beta_2 X' = 0.
\]
These equations suggest that a separated form is just as suitable for the other components of $\ten d$. Thus, suppose $u=u(x)$ and $v=v(y)$ to get
\[
\alpha_1 u'  + \alpha_2 X' = 0, \quad \beta_1 v' + \beta_2 Y' = 0, \quad
\alpha_2Y' +\beta_2 X' = 0.
\]
Accordingly,
\[
X' = \lambda \alpha_2, \quad Y' = -\lambda\beta_2,\quad
u' = -\lambda\frac{\alpha^2_2}{\alpha_1}, \quad v'=\lambda\frac{\beta^2_2}{\beta_1},
\]
up to some multiplicative factor (by linearity). Then, $(u,v,w=X+Y)$ can be obtained by direct integration:
\[
u = -\int\frac{\alpha^2_2}{\alpha_1}, \quad
v =  \int\frac{\beta^2_2}{\beta_1}, \quad
w =  \int\alpha_2 -\int\beta_2.
\]
Therein, $\alpha_1\neq0$ and $\beta_1\neq 0$ are guaranteed by hypothesis (2) stating that tangents to $\gt\alpha$ and $\gt\beta$ are never in each other's planes.

For the adopted examples, namely the ones in equation~\eqref{eq:ExSurf}, the deflection is best written in the canonical basis of $\R^3$. It reads
\[
\ten d = \begin{pmatrix}
    -\int\alpha^2_2 \\
    c\int\frac{\beta^2_2}{\beta_1}-s\left(\int\alpha_2 -\int\beta_2\right)\\
    s\int\frac{\beta^2_2}{\beta_1}+c\left(\int\alpha_2 -\int\beta_2\right)
\end{pmatrix},
\]
and is depicted in Figure~\ref{fig:deflections}. Here too, the case where the surface is periodic is particularly insightful as the stretching solution can be re-written as
\[
\ten d = \begin{pmatrix}
    -\average{\alpha^2_2}x \\
    \left(c\average{\beta^2_2/\beta_1}+s\average{\beta_2}\right)y\\
    0
\end{pmatrix} + \text{ periodic \& rigid motions}
\]
where $\average{\cdot}$ denotes taking the average of a periodic function over a period. Hence, on average, the effective stretching solution has an effective strain tensor
\[\label{eq:E(theta)}
    \ten E(\theta) = \begin{bmatrix}
        -\average{\alpha_2^2} & 0 \\ 0 & c\average{\beta^2_2/\beta_1}+s\average{\beta_2}
        \end{bmatrix},
\]
that is function of the inclination $\theta$. Strain $\ten E$ describes an effective biaxial stretch/contraction that is locally isometric: the strain in the $x$-direction is clearly a contraction; it is not clear if the strain in the $y$-direction is positive or negative. For $\theta=0$, e.g., for the eggbox pattern, the effective strain
\[
    \ten E(\theta) = \begin{bmatrix}
        -\average{f'^2} & 0 \\ 0 & \average{g'^2}
        \end{bmatrix},
\]
describes a stretch coupled to a contraction, a behavior characteristic of a positive Poisson's coefficient; see Figure~\ref{fig:deflections}a, b. However, for $\theta=\pi/2$, e.g., for the Miura ori, the effective strain
\[
\ten E(\theta) = \begin{bmatrix}
        -\average{f'^2} & 0 \\ 0 & -1
        \end{bmatrix},
        \]
describes a biaxial contraction, characteristic of an auxetic behavior, i.e., a negative Poisson's coefficient; see Figure~\ref{fig:deflections}c, d.

In general, as the plane of $\gt\alpha$ is rotated, there comes a point where it becomes tangent to $\gt\beta$. For that critical value of $\theta$, the denominator $\beta_1\equiv c+sg'$ vanishes and the the stretch in the $y$-direction changes signs thus switching between positive and negative Poisson's coefficient. This transition has been observed in the ``morph'' pattern, i.e., in cases where $f$ and $g$ are zigzags~\cite{Pratapa}. The phenomenon is quite general as is seen here and is encapsulated in the closed-form expression of the effective Poisson's coefficient
\[
    \nu(\theta) \equiv -\frac{E_{11}(\theta)}{E_{22}(\theta)} = \frac{\average{\alpha_2^2}}{c\average{\beta_2^2/\beta_1}+s\average{\beta_2}}=\begin{cases}
        \average{f'^2}/\average{g'^2}\quad&\text{for} \quad \theta=0,\\
        -\average{f'^2}\quad&\text{for}\quad \theta=\pi/2.
    \end{cases}
\]

\subsection{Bending}
The ansatz for bending $w(x,y) = X(x)Y(y)+S(x)+T(y)$ is more sophisticated: it combines elements of the stretching and twisting solutions. On one hand, it turns out that it is hard to anticipate whether a bending solution describes an in-plane or an out-of-plane state; when bending occurs in plane, it amounts to a gradient of stretching which is ultimately why a ``hint of stretching'' is needed in the ansatz. On the other hand, certain bent and twisted states are equivalent up to a change in basis which is why it is not surprising that the ansatz has a ``hint of twisting'' as well. In any case, substituting into the first two isometric constraints of~\eqref{eq:isoOfTrans} leads to
\[
    \alpha_1 u_x + \alpha_2(X'Y+S') = 0,\quad 
    \beta_1 v_y + \beta_2(XY'+T') = 0.
\]
These can be integrated as before into
\[
  u = -Y\int\frac{\alpha_2}{\alpha_1}X' - \int\frac{\alpha_2}{\alpha_1}S'-u_o,\quad 
  v = -X\int\frac{\beta_2}{\beta_1}Y' - \int\frac{\beta_2}{\beta_1}T'-v_o,
\]
where $u_o=u_o(y)$ and $v_o=v_o(y)$ are two integration ``constants''. The third isometric constraint of~\eqref{eq:isoOfTrans} then yields
\begin{multline}
- \alpha_1 \left(Y'\int\frac{\alpha_2}{\alpha_1}X'+u_o'\right) - \beta_1 \left(X'\int\frac{\beta_2}{\beta_1}Y'+v_o'\right) \\
+ \alpha_2 (XY'+T')+\beta_2(X'Y+S') = 0.    
\end{multline}
Particular choices are adopted here in an attempt to balance $x$-dependent and $y$-dependent terms. Recall that we are looking for particular solutions that correspond to bending rather than for the most general solution to these equations; see the discussion below. For our purposes then, four choices appear to be suitable and lead to four solutions:
\begin{itemize}
    \item Solution ($s$) is such that $X'=\alpha_2$, $Y'=\beta_2$;
    \item Solution ($p$) is such that $X'=\alpha_2$, $Y'=\beta_1$ and $S=0$;
    \item Solution ($\bar{p}$) is such that $X'=\alpha_1$, $Y'=\beta_2$ and $T=0$;
    \item Solution ($t$) is such that $X'=\alpha_1$, $ Y'=\beta_1$, $S=0$ and $T=0$.
\end{itemize}
Note that Solution~($t$) corresponds to twisting and has been explored above. As for the other three, the first one denoted $(s)$ is symmetric in its treatment of $\alpha$ and $\beta$ whereas the remaining two $(p)$ and $(\bar{p})$ permute the way in which $\alpha$ and $\beta$ appear. Upon integration, these solutions are found and listed hereafter:

\begin{itemize}
    \item Solution~$(s)$ is
\[\label{eq:sols}
    \begin{split}
    u &= -\int\beta_2\int\frac{\alpha_2^2}{\alpha_1} - \int\left(\alpha_2\int\frac{\alpha_2^2}{\alpha_1}-\frac{\alpha_2^2}{\alpha_1}\int\alpha_2\right),\\
    v &= -\int\alpha_2\int\frac{\beta_2^2}{\beta_1}-\int\left(\beta_2\int\frac{\beta_2^2}{\beta_1}-\frac{\beta_2^2}{\beta_1}\int\beta_2\right),\\
    w &= \int\alpha_2\int\beta_2+\int\left(\beta_1\int\frac{\beta_2^2}{\beta_1}-\beta_2\int\beta_2\right)+\int\left(\alpha_1\int\frac{\alpha_2^2}{\alpha_1}-\alpha_2\int\alpha_2\right).
    \end{split}    
\]
\item
Solution~$(p)$ is
\[\label{eq:solp}
    \begin{split}
    u &= -\int\beta_1\int\frac{\alpha^2_2}{\alpha_1},\\
    v &= -\int\alpha_2\int\beta_2-\int\left(\beta_2\int\beta_2-\frac{\beta_2^2}{\beta_1}\int\beta_1\right)-\int\left(\alpha_1\int\frac{\alpha_2^2}{\alpha_1}-\alpha_2\int\alpha_2\right),\\
    w &= \int\alpha_2\int\beta_1+\int\left(\beta_1\int\beta_2-\beta_2\int\beta_1\right).
    \end{split}    
\]
\item
Solution $(\bar{p})$ is obtained from solution $(p)$ by permuting $u$ and $v$ as well as $\alpha$ and $\beta$.
\end{itemize}

Two combinations of the solutions $(s)$, $(p)$ and $(\bar p)$ are depicted in Figure~\ref{fig:deflections}, one corresponding to out-of-plane bending, the other to in-plane bending. Formally, the fact that these solutions describe bending is evidenced by the presence of double integrations which produce quadratic terms in $x$ and $y$. But a finer analysis is needed to ascertain how bending takes place. To make this more precise, suppose $\gt\alpha$ and $\gt\beta$ are periodic oscillating about some average directions $\average{\gt\alpha'}$ and $\average{\gt\beta'}$. Then, the out-of-plane component is taken in direction $\ten N\equiv\average{\gt\alpha}\wedge\average{\gt\beta}$ and reads
\[
    \ten d\cdot\ten N = -\average{\alpha_2}\average{\beta_1}u-\average{\alpha_1}\average{\beta_2}v+\average{\alpha_1}\average{\beta_1}w.
\]
The out-of-plane bending solution is characterized by an out-of-plane component $\ten d\cdot\ten N$ that is quadratic, up to a periodic correction, in the cartesian coordinates $(\bar{x},\bar{y})$. Equivalently, the sought solution is free of secular terms that are products of periodic and linear or quadratic terms. For the eggbox pattern $(\theta=0)$, it is solution $(s)$ that describes out-of-plane bending (Figure~\ref{fig:deflections}a, b). For the Miura-ori $(\theta=\pi/2)$, it is solution $(p)$ (Figure~\ref{fig:deflections}c, d). In general, the out-of-plane bending solution is a linear combination of solutions $(s)$, $(p)$ and $(\bar p)$. For the adopted examples, given that $\average{\alpha_2}=0$, $\alpha_1=1$, and $\average{\beta_1}\beta_1+\average{\beta_2}\beta_2=1$, it appears that the linear combination with weights $(-\average{\beta_1},\average{\beta_2},0)$ of solutions $(s)$, $(p)$ and $(\bar p)$, respectively, is suitable in that it is free of secular terms. Upon simplification, the out-of-plane component reads
\[
\ten d\cdot\ten N = 
-\left(c\average{\beta_2^2/\beta_1} + s\average{\beta_2}\right)\bar{y}^2/2 -\average{\alpha_2^2}\bar{x}^2/2 \,\, +\text{ periodic \& rigid motions}.
\]
Accordingly, the \emph{quadratic behavior of the bending solution can be deduced from the linear behavior of the stretching solution}. In particular,
\[
\ten d\cdot\ten N \propto -E_{22}\bar{y}^2/2+E_{11}\bar{x}^2/2 \,\, + \text{ periodic \& rigid motions},
\]
where $(E_{11},E_{22})$ are the effective stretchings found in the stretching solution; see~equation~\eqref{eq:E(theta)}. In other words, \emph{the ratio of normal curvatures in the $x$- and $y$-directions in the bending solution yields the in-plane effective Poisson's coefficient in the stretching solution}, namely
\[
\frac{\average{\partial^2\ten d\cdot\ten N/\partial x^2}}{\average{\partial^2\ten d\cdot\ten N/\partial y^2}} = -\frac{E_{11}}{E_{22}} = \nu.
\]
Connections of this kind, between effective stretching and effective bending in isometric deformations of generally periodic surfaces, have been the subject of a recently developed theory~\cite{Nassar2025}. \emph{The crucial novelty here is in providing a closed-form expression for the periodic corrections, e.g., equations~\eqref{eq:sols} and \eqref{eq:solp}}. Such corrections had only been computed algorithmically in very simple cases, e.g., for the Miura-ori and other polyhedral origami tessellations with four parallelogram facets per unit cell. The present work computes these corrections for a vast, even if special, family of surfaces including polyhedral, smooth, curved-crease, developable, and doubly-curved surfaces. It therefore provides a simple and useful benchmark for the assessment of theoretical and numerical mechanical models that take into account elastic effects.

\section{Discussion}
We have shown that surfaces of translation admit three isometric deformations vaguely corresponding to twisting, stretching and bending. In the case where the surface is periodic, the solutions were rigorously shown to describe the twisting, stretching and bending kinematics of an equivalent Kirchhoff-Love plate, up to some periodic corrections free of secular terms. A recent theorem~\cite{Nassar2025} shows that periodic surfaces that are simply connected (i.e., with no holes or handles) admit exactly three solutions within these categories, i.e, with combined in-plane stretching and out-of-plane bending/twisting. Hence, the provided family of solutions is \emph{complete} in that regard.

Periodic surfaces typically admit other isometric deformations that do not fit into the kinematics of an equivalent plate. Examples mainly include:
\begin{itemize}
    \item Solutions with secular terms: these are solutions that contain products of linear terms (or unbounded terms more generally) and periodic terms. In other words, these are locally periodic solutions with a globally graded amplitude. In the present context, solutions $(s)$, $(p)$ and $(\bar p)$ span a 3D space of solutions of which only one dimension corresponds to out-of-plane bending. The other two dimensions are spanned by solutions with secular terms and correspond to ``bending in the plane'' or, equivalently and more appropriately, to a gradient of in-plane stretching. Such in-plane states are exemplified in Figure~\ref{fig:deflections}.
    \item Periodic solutions, i.e., with no linear or quadratic components. These solutions describe local isometric motions that leave no trace over a unit cell. They can be periodic with the same periodicity as the surface or with a larger period than the unit cell, i.e., with a doubling, tripling and so on of the unit cell. These solutions solve the linear equations of isometric deformations under generalized periodicity conditions of the Bloch type. In the simplest cases, the problem reduces to solving a pair of Hill's or even Mathieu's equations. Although these solutions leave no visible trace globally, they play a crucial role in the response of an elastic surface as they define ways in which the surface can relax under prescribed load conditions. The simplest example of such solutions, for polyhedral surfaces, is any deflection that is locally normal to facets and compactly supported over facets interiors.
\end{itemize}

\begin{sidewaysfigure}[!ht]
     \centering
     \includegraphics[width=\linewidth]{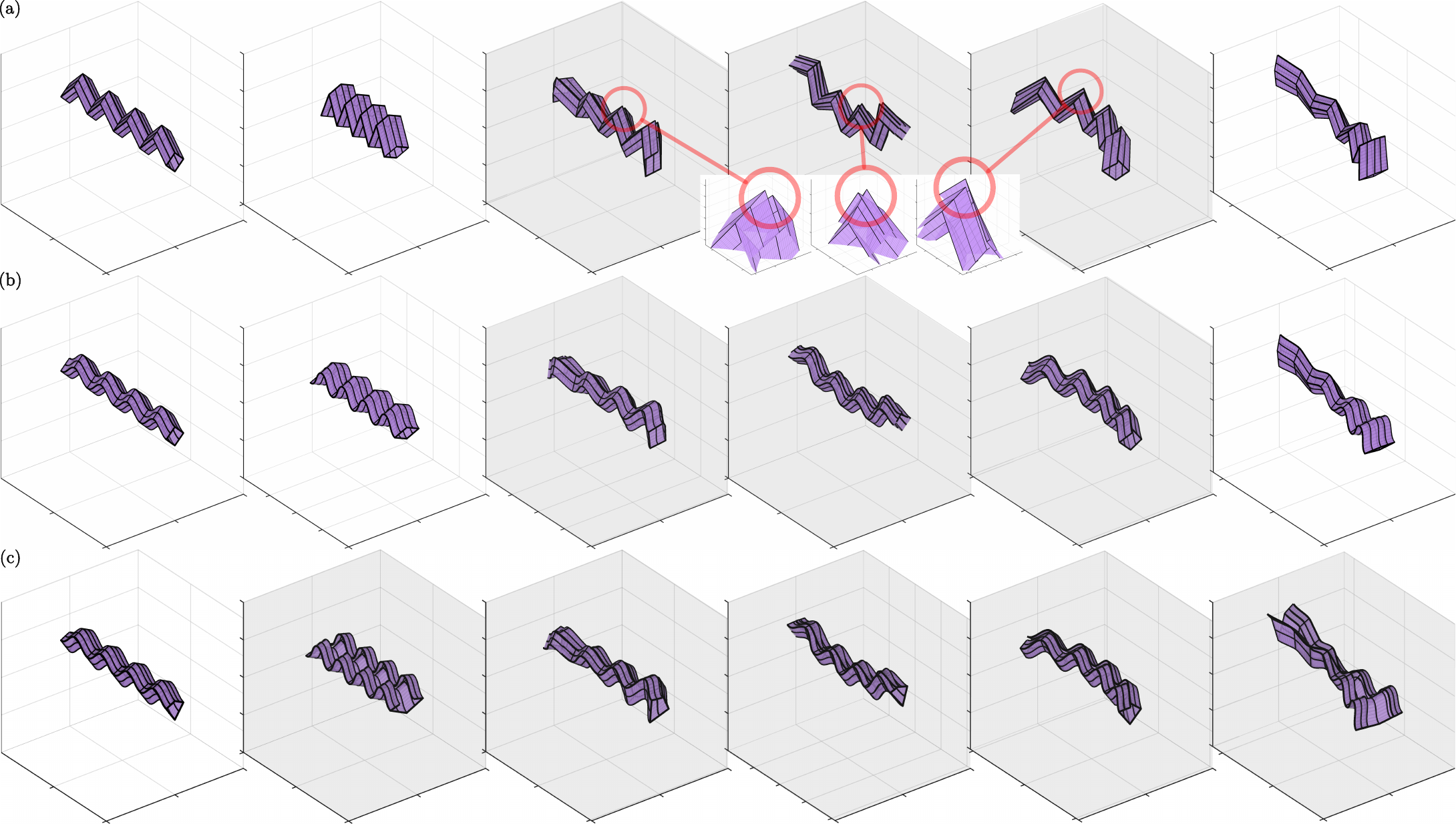}
     \caption{Infinitesimal isometries illustrated for tubular surfaces: (a) square cross-section with a zigzag path; (b) square cross-section with sinusoidal path; (c) non-centrosymmetric cross-section with sinusoidal path. Modes, per column, left to right are: reference, axial compression, torsion, bending in the vertical plane, bending in the horizontal plane, and axial stretch gradient. Close ups (red circles) reveal the presence of dislocations; all states with a grayed background have dislocations. The Matlab code that produced these figures is available at \url{github.com/nassarh}.}
    \label{fig:deflectionsTubular}
\end{sidewaysfigure}

It is crucial to recall that the derived solutions are not restricted by periodicity conditions. As an example, consider surfaces that describe tubular structures useful for the design of deployable booms. The twisting solution, still valid, now corresponds to torsion; stretching remains as is, with the Poisson effect now describing a change in the cross-section; and bending solutions $(s)$, $(p)$ and $(\bar p)$ now describe two states bent about the two axes of the cross-section and one state that is gradually stretched along the span of the tube; see Figure~\ref{fig:deflectionsTubular}. Where solutions might fail is when a new topology is enforced, e.g., that of a cylinder. For instance, it is well-known that Saint-Venant's solution to the problem of torsion for thin-walled prismatic bars involves stretching if the cross-section is closed and is isometric if the cross-section is open. Indeed, the topology of the latter case allows for the presence of a dislocation where the former precludes it. The dislocation here is the Saint-Venant warping integrated over the cross-section, i.e., the warping jump at the ``seam''. Therefore, all twisting solutions, for closed cross-sections, feature a discontinuity; see Figure~\ref{fig:deflectionsTubular} (column 3). In the same fashion, it appears that the bending solutions feature discontinuities as well (columns 4 and 5 of the same figure). For stretching, and stretching gradients, the appearance of dislocations depends on the cross-section (columns~2 and 6). In particular, the centro-symmetric cross-sections (a) and (b) have no dislocations in the stretched configuration but cross-section (c) does and lacks centro-symmetry.

\section{Conclusion}
Infinitesimal isometries inform on the stiffness of compliant shell mechanisms as well as on their local kinematics, even under large transformations. Here, for surfaces of translation, a \emph{complete} family of infinitesimal isometries has been derived in closed form by using an inspired scheme of variable separation and solving the local isometric constraints. The completeness refers to the fact that the family spans three out of the six deformation modes of an equivalent Kirchhoff-Love plate including stretching, bending and twisting; the other three modes have been proven to be non-isometric elsewhere~\cite{Nassar2025}. Solutions corresponding to stretching and twisting are relatively straightforward to obtain and have been known, in one form or another, see e.g.\ \cite{Nassar2023}. The bending solutions are far more challenging; their derivation in closed form is the first result of its kind. The construction reported here confirms, in a brute force fashion, certain properties of compliant shells, such as the equality of the effective Poisson coefficient to the ratio of normal curvatures under bending. Beyond effective properties, the construction further provides the underlying periodic correctors and thus allows to reconstruct candidate local fields from global ones.

\clearpage
\clearpage
\bibliographystyle{amsplain}
\bibliography{lib}
\end{document}